\documentclass[11pt]{article}
\usepackage{amsmath,amsthm,amsfonts,amssymb,mathrsfs,epsfig,bm}
\numberwithin{equation}{section}

\parskip        \medskipamount
\newtheorem{theorem}{THEOREM}[section]
\newtheorem{lemma}{LEMMA}
\newtheorem{corollary}{Corollary}

\def\Z{\mathbb{Z}}
\def\R{\mathbb{R}}
\def\D{\mathbb{D}}

\def\P{\mathbb{P}}
\def\E{\mathbb{E}}
\def\O{\mathbb{O}}

\renewcommand{\phi}{\varphi}
\renewcommand{\epsilon}{\varepsilon}
\newcommand{\cov}{\operatorname{cov}}

\renewcommand{\liminf}{\varliminf}

\long\def\symbolfootnote[#1]#2{\begingroup
\def\thefootnote{\fnsymbol{footnote}}\footnote[#1]{#2}\endgroup}
\newcommand{\keywords}[1]{ \noindent {\footnotesize
             {\small \em Keywords and phrases.} {\sc #1} } }
\newcommand{\ams}[2]{  \noindent {\footnotesize
             {\small \em AMS {\rm 2010} subject classification.
             {\rm Primary {\sc #1}; secondary {\sc #2}} } } }

\def\Spin{\mathsf S}

\def\I{\textbf{I}}

\theoremstyle{plain}
\newtheorem*{lemma*}{Lemma}
\newtheorem*{theorem*}{Theorem}

\def\dd{\bm d}
\def\ent{\operatorname{\mathsf{Ent}}}

\newcommand{\III}[1]{{\left\vert\kern-0.5ex\left\vert\kern-0.5ex\left\vert #1 
    \right\vert\kern-0.5ex\right\vert\kern-0.5ex\right\vert}}

\begin{document}

\title{\Large Log-Sobolev inequalities for infinite-dimensional Gibbs measures with non-quadratic interactions} 
\author{James Inglis$^*$
 \and  Ioannis Papageorgiou$^{**}$}

\date{ }

\maketitle

\begin{abstract}
We focus on the log-Sobolev   inequality for spin systems on the lattice with   interactions of  higher order than quadratic.
We show that if the one-dimensional single-site measure with boundaries
satisfies the log-Sobolev inequality uniformly in the boundary conditions
then the infinite-dimensional Gibbs measure   also satisfies  the inequality
under   appropriate conditions on the  phase and the interactions.

\symbolfootnote[0]{
\textit{Address:}  $^{**}$ Neuromat, Instituto de Matematica e Estatistica,
 Universidade de Sao Paulo, 
 rua do Matao 1010,
 Cidade Universitaria, 
 Sao Paulo,  SP  05508-090,  Brasil.
\\ \text{\  \   \      } \textit{Email:}$^*$inglis.jd@gmail.com $^{**}$  ipapageo@ime.usp.br, papyannis@yahoo.com   \\
 $^{**}$
This article  is supported by FAPESP grant (2017/15587-8);  
This article was produced as part of the activities of FAPESP  Research, Innovation and Dissemination Center for Neuromathematics (grant $\#$2013/ 07699-0 , S.Paulo Research Foundation).}

\vspace*{2mm}

\vspace*{2mm}
\ams{60K35,   
39B62,          
26D10, 
}
{
82B20           
}

\end{abstract}

%%%%%%%%%%%%%%%%%%%%%%% BEGIN SECTION INTRODUCTION %%%%%%%%%%%%%%%%%%%%%%%%%
\section{Introduction}
We focus   on the    log-Sobolev   
inequality  for probability measures on unbounded spin systems on the lattice with  nearest neighbour interactions of power which is  higher than quadratic.  
Our  objective  is to determine conditions on the phase and interactions   
so that the log-Sobolev inequality can be extended from the one site  measure with a boundary
to the infinite dimensional Gibbs measure.
The main assumption of the paper  is that the single-site log-Sobolev inequality holds  
with a constant which is independent of  the boundary conditions, while  we assume that the phase is of higher power than the  interaction.

We will denote $\Spin$  the single-site space and the state space $\Omega:={\Spin}^{\Z^d}$, while with $\P^{\Lambda,\omega}$ the probability measure on $\Spin^\Lambda$, 
for any finite  $\Lambda \subset \Z^d$. These are   the  local specifications, they
 depend on the boundary conditions $\omega \in \Spin^{\partial \Lambda}$ and 
satisfy the usual spatial Markov property.  They are of Gibbs type  and their   Hamiltonian 
is composed by  two parts: the phase that depends on a single
sites and the interaction that depends on two neighbouring sites.
The  integration with respect 
to $\P^{\Lambda,\omega}$ is denoted by $\E^{\Lambda,\omega}$.  
The log-Sobolev inequality for  local specifications with quadratic interactions
  has been studied by \cite{G-Z},  \cite{Led},
  in \cite{Z1} and \cite{Z2}, 
  \cite{B-E}, 
  \cite{A-B-C},
 \cite{B-H} and  \cite{Y}.
Furthermore, in  \cite{G-R} the Poincare inequality 
is studied. 
For the one node measure  
$\mathbb{E}^{\{i\},\omega}$, 
 conditions for the logarithmic Sobolev inequality are presented in  \cite{R-Z},  \cite{B-G}, 
 \cite{B-Z}.
Furthermore,  the question  of passing from the uniform inequality of the single-site measure $\mathbb{E}^{\{i\},\omega}$  directly to the inequality about the  infinite-dimensional Gibbs measure
when the interactions are  quadratic has been 
 studied  by  \cite{M},  
 \cite{O-R} and  \cite{I-P}.

In the current paper we are interested in   the case of    interactions that grow faster than a quadratic. 
The case of higher than  quadratic interactions has been investigated before in \cite{Pa1}, but only for spin systems on the one-dimensional lattice.
Now we determine conditions so that the inequality can be extended from the one site measure to the Gibbs measure on the $d$ dimensional lattice.

%%%%%%%%%%%%%%%%%%%%%%%%%%%%%%%%%%%%%%%%%%%%%%
\subsection{General framework}
We consider a $d$ dimensional integer lattice $\Z^d$ which is equipped with 
the typical neighborhood
structure, where  for two  neighbouring    sites, say   $i, j \in \Z^d$,  we will write $i \sim j$. If  $\Spin$ is the  spin space,  we will work with   the configuration
space $\Omega = \Spin^{\Z^d}$, so that  the coordinate $x_i$ of a configuration $ x \in \Omega$ will correspond  to the spin at site $i$, with $x_i$ taking values in $\Spin^i \equiv \Spin$.
For every subset $\Lambda \subset \Z^d$ we  then identify $\Spin^\Lambda$
with the Cartesian product of the  $\Spin^i$, with $i$ ranging over $\Lambda$.
Furthermore, we assume that the spin space  $\Spin^{i}$ has a natural measure which we denote as $dx_i$,  while for the  product measure  of the $dx_i$, $i \in \Lambda$ we will  write $dx_\Lambda$. As an  example,
in the case of a group   $\Spin$,  one can  think  $dx_i$ as the invariant under
the group operation measure.
  Concerning the measure on the spin space  $\E^{\{i\},\omega}$ we assume that it is absolutely continuous
with respect to $dx_i$. Then, because of the  Markov property,
for any finite subsets $\Lambda$ of $\Z^d$, there exists a Hamiltonian $H^{\Lambda, \omega}$ (see \cite{Pr}) so that the probability measures
$\E^{\Lambda, \omega}$  have the following   form:
\[
\E^{\Lambda, \omega} (dx_\Lambda) = \frac{e^{-H^{\Lambda,\omega}(x_\Lambda)}\, dx_\Lambda}{Z^{\Lambda,\omega}}\,
,
\]
where $Z^{\Lambda,\omega}$ is the normalization constant, 
while  the Hamiltonian $H^{\Lambda, \omega}$  is of 
the form  
\[
H^{\Lambda,\omega}(x_\Lambda) 
:= \sum_{i \in \Lambda} \phi(x_i) 
+ \sum_{i,j \in \Lambda,\, j \sim i} J_{ij} V(x_i, x_j)
+ \sum_{i\in \Lambda,j \in \partial \Lambda,\, j \sim i} J_{ij} V(x_i, \omega_j),
\]   
 For a function $f:\Spin^{\Z^d}\rightarrow \R$,
we conventionally write $\E^{\Lambda, \omega} f$ for  the expectation of $f$ with respect to the measure  $\E^{\Lambda, \omega} $, obtained by integrating with respect
to  $dx_\Lambda$ while  substituting the neighbouring nodes $x_{\partial \Lambda}$ by the boundary conditions  $\omega$. 
For simplicity we will frequently write   $\E^\Lambda f$ instead of $\E^{\Lambda, \omega} f$.
 The  Markov property  takes the following expression
\[
\E^\Lambda \E^K = \E^\Lambda, \quad K \subset \Lambda.
\]

We define the  infinite volume Gibbs measure $\nu$   on $\Omega = \Spin^{\Z^d}$
 as the probability  for the local
specifications $\{\E^{\Lambda,\omega}\}$ which satisfies 
the Dobrushin-Lanford-Ruelle equations:
\[
\nu \mathbb{E}^{\Lambda,\bullet}=\nu, \quad \Lambda \Subset \Z^d,
\]
which means that $\nu$ is an invariant measure for the Markov
random field (see 
  \cite{Pr},  \cite{B-HK} and   \cite{D}).
Throughout the paper we   assume both the existence and uniqueness of   $\nu$, even though  uniqueness is 
deduced from our main results.

We now present the main framework   about 
the spin space $\Spin$. 
We consider a     spin space $\Spin$ that is  a nilpotent  Lie group on
$\R^d$ with a H\"ormander system $X^1, \ldots, X^n$, $n \le d$, for which we assume that 
if the vector fields are $X^k = \sum_{j=1}^d a_{kj} \frac{\partial}{\partial x_j}$
$k=1,\ldots, n$, then the coefficients  $a_{kj}$ are functions of $x \in \R^d$
that do not  depend on the $j$-th coordinate $x_j$. 
We then define the (sub)gradient $\nabla$ with respect to 
this system of vector fields to be   $\nabla f = (X^1 f, \ldots, X^n f)$,
and the (sub)Laplacian $\Delta = (X^1)^2 + \cdots + (X^n)^2$. 
Then  $\|\nabla f\|^2 := (X^1 f)^2 + \cdots + (X^n f)^2$.
We denote  $\nabla_i$ and $\Delta_i$ the gradient and Laplacian respectively acting on  functions on the spin space $\Spin^i, i \in \Z^d$. For a finite 
 $\Lambda \subset \Z^d$ we then define 
$\nabla_\Lambda := (\nabla_i, i \in \Lambda)$ and 
$\|\nabla_\Lambda f\|^2 := \sum_{i \in \Lambda} \|\nabla_i f\|^2$.
Furthermore, it is assumed  that $\Spin$ is equipped with a metric-like
function $\dd(x,y)$, for $x,y \in \Spin$.
For instance, in the case where  $\Spin$ is a Euclidean space,  $\dd$ is the Euclidean
metric, while in the case where  $\Spin$ is the Heisenberg group, 
 $\dd$ is the Carnot-Carath\'eodory metric. 
We adopt the following convention; for $x \in \Spin$, we  shall  denote  the metric \[\dd(x):=\dd(x, 0),\] for some specific point  $0$  of $\Spin$, for
example  in the case where    $\Spin$ is $\R^m$, $0$ can be the origin while, if $\Spin$ is a Lie group  the identity element of $\Spin$.

The main hypothesis of the paper is that the probability measures 
$\E^{\{i\},\omega}$ 
relate to the differential structure
via a log-Sobolev inequality on $\Spin^i$ with a
constant which does not depend on $\omega$, that is,
that there exists $c>0$ such that
\begin{equation}
\label{ei}
\E^{\{i\},\omega}\bigg( f^2 \log \frac{f^2}{\mathbb{E}^{\{i\},\omega} f^2}
\bigg) \le c\, \mathbb{E}^{\{i\},\omega} \|\nabla_i f\|^2,
\quad i \in \Z^d, \, \omega \in \Spin^{\partial \{i\}},
\end{equation}
for any smooth function $f: \Spin^i \mathsf \to \R$.

We say that the log-Sobolev inequality holds for
$\E^{\{i\},\omega}$ {\em uniformly} (in $\omega$.)
We notice that assumption \eqref{ei} does not imply 
that the log-Sobolev inequality holds with the same constant $c$
for the measures $\E^{\Lambda,\omega}$ even when
$\Lambda$ is a finite subset of $\Z^d$. 
If, however, $\Lambda$ is a subset  of $\Z^d$ such that
any two points of $\Lambda$ have\ distance   greater or equal to two from each other, then the log-Sobolev inequality holds for
$\E^{\Lambda,\omega}$, with the same constant $c$, uniformly
in $\omega \in \partial \Lambda$.

Similar thing is also true for 
spectral gap inequalities (a measure $\mu$ satisfies spectral gap inequality
with constant $C$ if
$\mu |f-\mu f|^2 \le C\, \mu |\nabla f|^2$). For proofs of these properties one   can look  in  
   \cite{G-Z},   \cite{G} and 
 \cite{B-Z}.

\section{Main assumptions and  results }

This section concerns presenting the main assumptions on the phase and interactions and the principle results about the infinite volume Gibbs measure.

 We recall, that the one site Hamiltonian is defined by 
\[
H^{i}(x_i) 
:=  \phi(x_i) 
+ \sum_{ j \sim i} J_{ij} V(x_i, \omega_j),
\]
where above for economy we suppressed the boundary $\omega$, writing $H^i$ for $H^{i,\omega}(x_i) $. Of course, we assume that the  $\phi$ and $V$
are such that $\int_\Spin \exp(-H(x_i)) dx_i < \infty$. 

\subsubsection*{The main assumption}$~$
The  one-dimensional measures $\E^{i,\omega}$ satisfy
the log-Sobolev inequality with a constant $c$ uniformly with respect to the
boundary
conditions $\omega$.

\subsubsection*{Assumptions on the  local specification}
We  assume that $\phi$ and the $V$ are twice  continuously
 differentiable and that there exist
  nonnegative 
$C^2$ functions  $\hat \phi$  and $  \hat V(x_i,x_j)$
such that
\begin{equation}
\label{Hphi}
\nabla_i \phi (x_i) = \hat \phi (x_i)\nabla_i \dd (x_i)
\end{equation}
and
\begin{equation}
\label{HV}
\nabla_i V(x_i,x_j)=   \hat V(x_i,x_j) \nabla_i \dd(x_i).\end{equation}
 There exist constants $\xi$ and $\zeta$ such that, for all
$x_i \in \Spin^i$,
\begin{equation}
\label{Hgrad}
\zeta \le \|\nabla_i \dd (x_i)\| \leq \xi
\end{equation}
and a constant $\beta$ with
\begin{equation}
\label{Hlap}
|\Delta_i \dd (x_i) | \le \frac{\beta}{\dd (x_i)}.
\end{equation}
Moreover, we require that there exists $k_1 >0$ and 
$ q \ge 2$ 
such that
\begin{equation}
\label{Hmetric}
   \dd^q(x_i) \leq \phi (x_i)   \ \text{and} \   k_1 \dd^q(x_i)\le \dd(x_i)\hat \phi (x_i)  .
\end{equation}

Furthermore, the interaction  potential $V$ depends on $x_i, x_j \in \mathsf \Spin$
only through $\dd(x_i)$ and $\dd(x_j)$ in a way that
 $\exists ~  r \leq q$ such that there exists a $\lambda>0$ 
   \begin{equation}
\label{HgradV2}
\E^je^{\epsilon \|\nabla_i V\|^2}
\le e^{\lambda +\lambda  \sum_{\ell\sim j} \dd(x_\ell)^r}
\end{equation}
for some $\epsilon >0$, and
\begin{equation}
\label{HgradV1}
\|\nabla_i V(x_i,x_j)\|^2 \le \lambda + \lambda \dd(x_i)^r +\lambda  \dd(x_j)^r
.\end{equation}
 We can now present the main theorem of the paper.
\begin{theorem}
\label{thmGENERAL} 
Let $f\colon \mathbb{\mathsf{M}}^{\mathbb{Z}^d} \to \mathbb{R}$. 
Assume that  (\ref{Hphi})-(\ref{HgradV2}) hold
and that  the one-dimensional measures $\mathbb E^{i,\omega}$ satisfy
a log-Sobolev inequality uniformly on the boundary conditions. 
 Then
$\nu$ satisfies a log-Sobolev inequality:
\[
\nu f^{2}\log\frac{f^{2}}{\nu
f^{2}}\leq \mathfrak{C} \ \nu \left\| \nabla f
\right\|^2,
\]
for some positive  constant $\mathfrak{C}$. 
\end{theorem} 
\addcontentsline{toc}{subsubsection}{The main theorem}

The main assumption about the one site  
measure $\mathbb{E}^{i,\omega}$ is that it satisfies the log-Sobolev
inequality with  a constant uniformly to the boundary conditions. Furthermore,  we require  the phase $\phi$ to  dominate over
the interactions, in the sense that
$$\left\| \nabla _{j}V(x_{i},\omega_{j})\right\|^{2}\leq \lambda +\lambda (
d^r(x_{i})+d^r(\omega_j))\leq \lambda +\lambda ( \phi(x_i
)+ \phi(\omega_j ))$$
for $r\leq q$.  As examples of Hamiltonians that satisfy conditions  (\ref{Hphi})-(\ref{HgradV2}), one can think $\mathsf{M}$ to be the Heisenberg group,  phase $\phi (x)=\dd^r(x)$ and interactions $V(x,y)=(\dd(x)+\dd(y))^s$ or $V(x,y)=\dd^s (x+y)$, with $s<r$.

\section{U-bound  inequalities.}
We prove a single-site coercive inequality, inspired by the U-bound inequalities introduced in \cite{H-Z}.  This will be an important tool in controlling the higher order interactions. 
 \begin{lemma}\label{coers}
   Assume (\ref{Hphi})-(\ref{HgradV2}). Then there exists a positive constant $G_{0}$ such that,
for all $r \le q$,
\[\E^i \dd^r f^2 \leq  G_{0} \E^i  |\nabla_i f|^{2}+G_{0} \E^i f^{2}\]
for any smooth function $f$ with compact support.
\end{lemma}
\begin{proof}
At first we notice that it is sufficient to show  the inequality for $r=2(q-1)$, since
 if $\E  \dd^{2(q-1)} f^2 \le K \E |\nabla f|^2 + K \E f^2$,
 for some $K>0$ then for all $r \le 2(q-1)$ we have
\begin{align*}\E^i \dd^r f^2
\leq
\E^i f^2 + \E^i \dd^{2(q-1)} f^2 \leq   K
\E^i |\nabla_i f |^2 + (K+1) \E^i f^2\end{align*}
If we   calculate 
\[
\nabla_i (f e^{-H^i}) = (\nabla_i f) e^{-H^i} 
-  (\nabla_i H^i) e^{-H^i} f
\]
and then 
 take on both sides the inner product  with $\dd\nabla_i \dd$,  we have
\[
 \dd\langle \nabla_i \dd,\nabla_i H^i\rangle e^{-H^i } f=\dd\langle\nabla_i \dd,\nabla_i f\rangle e^{-H^i}-\dd\langle \nabla_i \dd, \nabla_i (f e^{-H^i})\rangle 
\]
After  integrating with respect to $dx_i$ and dividing with the normalization constant $Z^i=\int e^{-H^i(x_i)} dx_i$ we obtain
\[
\I_{1}:=\E^i \dd\langle \nabla_i\dd,\nabla_i H^i\rangle f =\E^i\dd\langle\nabla_i\dd,\nabla_i f\rangle-\frac{1}{Z^i} \underbrace{\int \dd\langle \nabla_i\dd, \nabla_i (f e^{-H^i})\rangle dx_i}_{\I_2}\]
Since \[\E^i\dd\langle\nabla_i\dd,\nabla_i f\rangle \leq 
  \E^i\dd|\nabla_i\dd||\nabla_i f| \leq \xi
  \E^i\dd |\nabla_i f|,\] because of  (\ref{Hgrad}), 
  we get
  \[\I_{1}\leq  \xi
  \E^i\dd |\nabla_i f|  -  \frac{1}{Z^i}  \I_2.\]
Because of the structural assumptions about any of the H\"ormander generators of $\Spin^i$, for a generator  say $X^j$, we have  the following integration by parts
formula
\[
\int G (X^j F) dx_i = - \int (X^jG) F dx_i
\]
for smooth functions $F,G$ with  compact support.  As a consequence, we obtain the following 
  integration by parts formula
\[
\int f \langle \nabla_i \Phi, \nabla_i \Psi \rangle dx_i
= - \int \langle \nabla_i \Phi, \nabla_i f \rangle \Psi dx_i
- \int (\Delta_i\Phi) \Psi f dx_i
\]
from which we get
\begin{align*}
\I_2 
 =-\int \dd|\nabla_i\dd| f e^{-H^i} dx_i-\int \dd(\Delta_i\dd)f e^{-H^i} dx_i \geq- \xi \int \dd fe^{-H^i}dx_i-\beta \int fe^{-H^i}dx_i
\end{align*}
because of  (\ref{Hgrad}) and  (\ref{Hlap}). For the first  term $\I_1$ we have
\begin{align*}
\I_1 
& = \E^i\dd(x_i)\langle \nabla_i\dd(x_i),\nabla_i\phi(x_i)\rangle f+\sum_{j\sim i}J_{ij}\E^i\dd(x_i)\langle \nabla_i\dd(x_i),\nabla_i V(x_i,\omega_j)\rangle f= 
\\
&= \E^i \dd(x_i)\hat \phi(x_i)  \vert \nabla_i\dd(x_i)\vert f+\sum_{j\sim i}J_{ij}\E^i\dd(x_i) \hat V(x_i,\omega_j)\vert \nabla_i\dd(x_i)\vert f 
\end{align*}
because of  (\ref{Hphi})-(\ref{HV}). If we use  a (\ref{Hgrad}) and (\ref{Hmetric}) we get
\[ \I_1  \geq \zeta\ k_1\E^i\dd(x_i)^q f
\] Combining all the above,   we obtain
\[
\E^i( \dd(x_i)^qf)\leq\ \frac{1}{\zeta\ k_1}
  (\xi
  \E^i\dd |\nabla_i f|+ \xi \E^i\dd f+\beta \E^i f) ) .\]
 If in the last inequality  we replace $f$ by $f^2$  we get
\begin{align*}
\E \dd^q f^{2}\leq & 
  \frac{1}{\zeta k_1}(2\xi\E^i\dd f |\nabla_i f|+ \xi \E^i\dd f^{2}+\beta \E^i f^{2} )\leq\frac{1}{\zeta k_1}(\xi\E^i  |\nabla_i f|^{2}+\xi\E^i\dd^{2}f^{2}+\xi\E^i\dd f^{2}+\beta \E^i f^{2})\end{align*}
  from the   Cauchy-Swartz inequality. Then
  \begin{align*} \E \dd^q f^{2}\leq & \frac{1}{ \zeta k_1}\left(\xi\E^i  |\nabla_i f|^{2}+\xi\E^i(\mathrm{I}_{\{\frac{4\xi}{\zeta k_1}\leq\dd^{q-2}\}} +\mathrm{I}_{_{\{\frac{4\xi}{\zeta k_1}>\dd^{q-2}\}}})\dd^{2}f ^{2}\right )\\ &+\frac{1}{ \zeta k_1}\left( \xi\E^i(\mathrm{I}_{\{\frac{4\xi}{\zeta k_1}\leq \dd^{q-1}\}}+\mathrm{I}_{\{\frac{4\xi}{\zeta k_1}>\dd^{q-1}\}})\dd f^{2}+ \beta \E^i f^{2}\right ) \\  \leq&\frac{1}{ 2}\E^i (\dd^q f^2) +\frac{1}{\zeta k_1}\left(\xi\E^i  |\nabla_i f|^{2}+ \left(4\xi(\frac{4\xi}{\zeta k_1})^\frac{2}{q-2}+4\xi(\frac{4\xi}{\zeta k_1})^\frac{2}{q-1}+ \beta\right) \E f^{2}  \right)
.\end{align*}
 We finally obtain
\begin{align*}
\E\dd^{q}  f^{2}\leq \frac{2}{\zeta k_1}\left(\xi\E^i  |\nabla_i f|^{2}+ \left(4\xi(\frac{4\xi}{\zeta k_1})^\frac{2}{q-2}+4\xi(\frac{4\xi}{\zeta k_1})^\frac{2}{q-1}+ \beta\right) \E f^{2}  \right)
\end{align*}
which proves the  lemma.
\end{proof}

\section{Sweeping out inequalities.}

We prove the first set of sweeping out inequalities. The function $\E^i f$
is   a function on $\Z^d$ that depends  on the variables
$x_j$ with $j\sim i$, and the input of the function $f$, excluding of course $x_i$.
We  have
\begin{lemma}\label{sweep1}
Assume that (\ref{Hphi})-(\ref{HgradV2}) hold.  Let $j\sim i$ for $i,j \in \Z^d$.
Then there exist constants $Q_1>0$ and $0<Q_2 < 1$ such that
\[
\nu \| \nabla_j (\E^i f)\|^2 \le Q_1 \nu \| \nabla_j f\|^2
+ Q_2 \nu \| \nabla_i f\|^2.
\]
\end{lemma}
\begin{proof}
Fix $i ,j \in \Z^d$ neighbouring sites. If we consider   $X_j^1, \ldots, X_j^n$ to be the H\"ormander
system for $\Spin^{\{j\}}$, we have
 $\|\nabla_j (\E^i f)\|^2 = \sum_{\alpha=1}^n (X^\alpha_j f)^2 $.
If we denote  $\theta_i=e^{-H^i}/Z^i$  the 
density of $\E^i$ with respect to the measure $dx_i$, we have
\begin{align}
(X^\alpha_j (\E^i f))^2 
&=\left( 
\int \theta_i (X^\alpha_j f) dx_i
+ \int (X^\alpha_j \theta_i) f dx_i  
\right)^2
\leq 
2 \E^i (X^\alpha_j f)^2
%2\left(\int \theta_i (X^\alpha_j f) dx_i\right)^2 
+ 2\left( \int (X^\alpha_j \theta_i) f dx_i\right)^2,
\label{F1a}
\end{align}
where  in the first term we used Jensen's inequality. 
If we take the sum over $\alpha$ and then  integrate with respect to the Gibbs measure $\nu$, the first term
on the right gives
$$\nu\sum_{a=1}^n (\E^i X^\alpha_j f)^2= \nu \|\nabla_j f\|^2.$$ 
Concerning the   second term in the right hand side of (\ref{F1a}), one  notices that both  $H^i$ and $Z^i $ depend on  variables $x_r$ where   $r$ ranges over all the neighbors
of $i$, including $j$. So,  we can compute 
 \[X_j \theta_i = -\frac{e^{-H^i}}{Z^i} (X_j H^i)
+ \frac{e^{-H^i}}{Z^i} \int \frac{e^{-H^i}}{Z^i}  (X_j H^i) dx_i 
= -\theta_i (X_j H^i) + \theta_i \int \theta_i (X_j H^i) dx_i.\]
From this  and using Leibniz' rule, 
we  obtain
\begin{equation}
\label{F1b}
\int (X^\alpha_j \theta_i) f dx_i
= - \E^i[(f-\E^i f)\, (X^\alpha_j H^i)].
\end{equation} 
At this point, we can apply Jensen's inequality once more to pass the square inside the expectation 
\begin{equation}
\label{F1c}
\left( \int (X^\alpha_j \theta_i) f dx_i\right)^2
\le \E^i[(f-\E^i f)^2\, (X^\alpha_j H^i)^2].
\end{equation}
Therefore, 
\begin{align*}
\sum_{\alpha=1}^n \left( \int (X^\alpha_j \rho_i) f dx_i\right)^2
\le J_{ij}^2 \, \E^i \big[ (f-\E^i f)^2 \|\nabla_j V(x_i,x_j)\|^2\big]
.\end{align*}
If we use (\ref{HgradV1}) to bound the interaction
\begin{align*}
\sum_{\alpha=1}^n \left( \int (X^\alpha_j \rho_i) f dx_i\right)^2
\le \lambda J_{ij}^2 \E^i (f-\E^if)^2 
+  \lambda  J_{ij}^2 \E^i (f-\E^if)^2 \dd(x_i)^r
+ \lambda  J_{ij}^2 \E^i (f-\E^if)^2 \dd(x_j)^r
\end{align*}
Next, if we sum up the first display of this proof over $\alpha$
and integrate with respect to the Gibbs measure   $\nu$ we get 
\begin{multline}
\label{nubla}
\nu \|\nabla_j(\E^i f)\|^2 \le 2 \nu \|\nabla_j f\|^2
+ 2\lambda   J_{ij}^2 \nu[ (f-\E^if)^2 ]
\\
+  2\lambda   J_{ij}^2 \nu \E^i[ (f-\E^if)^2 \dd(x_i)^r]
+ 2\lambda   J_{ij}^2 \nu \E^j[ (f-\E^if)^2 \dd(x_j)^r].
\end{multline}where above we used $\nu \E^i=\nu \E^j =\nu$.
We can bound the last two terms from the single-site  coercive inequality of Lemma \ref{coers}. The first one becomes
\[
\nu \E^i [ (f-\E^if)^2 \dd(x_i)^r]
\le G_0\, \nu \|\nabla_i f\|^2 + G_0\, \nu [(f-\E^if)^2].
\]
Similarly, we compute the second one
\begin{multline}
 \nu \E^j [ (f-\E^if)^2 \dd(x_j)^r]
\le G_0\, \nu \|\nabla_j (f-\E^i f)\|^2 + G_0\, \nu [(f-\E^if)^2],
\\ 
\leq 2G_0\, \nu \|\nabla_j f\|^2
+ 2G_0\, \nu \|\nabla_j(\E^i f)\|^2
+ G_0\, \nu [(f-\E^if)^2].
\label{applyssc}
\end{multline}
So, substituting these   two into \eqref{nubla} gives
\begin{multline*} 
\nu \|\nabla_j(\E^i f)\|^2  
\le (2+4\lambda J_{ij}^2 G_0)\nu\|\nabla_j f\|^2+  2\lambda J_{ij}^2 (1+2G_0) \nu [(f-\E^if)^2] 
+2\lambda J_{ij}^2 G_0\, \nu\|\nabla_i f\|^2
\\
%+4\lambda J_{ij}^2 G_0\, \nu\|\nabla_j f\|^2
+ 4\lambda J_{ij}^2 G_0\, \nu \|\nabla_j(\E^i f)\|^2
\end{multline*}
while from  the spectral gap inequality we can bound the variance, 
$\nu [(f-\E^if)^2] = \nu \E^i [(f-\E^if)^2] \le 
c \nu \E^i \| \nabla_i f\|^2 = c \nu \| \nabla_i f\|^2$,
where $c$ the log-Sobolev constant of the main hypothesis, and so
\[
(1-4 \lambda G_0 J_{ij}^2) \, \nu \|\nabla_j(\E^i f)\|^2
\le (2+4\lambda J_{ij}^2 G_0)\nu\|\nabla_j f\|^2
+ 2\lambda J_{ij}^2 (c+2cG_0 + 2G_0) \nu \| \nabla_i f\|^2.
\]
For $J < (2\lambda (c+2cG_0+4G_0))^{-1/2}$,  the inequality holds with
$Q_1 := (2+4\lambda G_0 J^2)/(1-4\lambda G_0 J^2)$
and
$Q_2 := 2\lambda (c+2cG_0+2G_0)J^2/(1-4\lambda G_0 J^2)<1$.
\end{proof}
From this lemma the next corollary follows.
\begin{corollary}\label{coers2}
Suppose that  (\ref{Hphi})-(\ref{HgradV2}) hold. Then for   if $i, j$ are neighbors in $\Z^d$, there exists some  $Q_3 > 0$, so that
\[
\nu[(f-\E^i f)^2 \dd(x_j)^s] \le  Q_3 \nu \|\nabla_j f\|^2
+ Q_3 \nu\|\nabla_i f\|^2.
\]\end{corollary}
\begin{proof}
At first bound $\nu\|\nabla_j(\E^i f)\|^2$ in  the right-hand side of \eqref{applyssc} 
by its upper bound from the inequality in the statement of Lemma \ref{sweep1},
and then replace the last term by its upper bound $c \nu \| \nabla_i f\|^2$ 
from the spectral gap inequality.
\end{proof}

%and consider the operator
%\[
%R:= \E^{\Gamma_{d-1}} \cdots \E^{\Gamma_0}.
%\]
%The meaning of this operator should be clear from the
%convention \eqref{convention}.
In the remaining of the paper   we shall work with the two dimensional lattice for simplicity of notation.
The general case is analogous. Next, let, for $r =0,1$, the set $\Gamma_r$
 defined by
\[
\Gamma_r := \{i \in \Z^2:\, i_1+i_2 \equiv r \mod 2\}.
\]
Note that the sets $\Gamma_r$, $r =0,1$,
form a partition of $\Z^2$ and
$\inf\{\max_{1\le k\le 2} |i_k-j_k|:\, i \in \Gamma_r, j \in \Gamma_s\}=1$
if $r \neq s$.

\begin{lemma}\label{sweepGam1}Suppose that  (\ref{Hphi})-(\ref{HgradV2}) hold.
%Let $\Gamma_0, \Gamma_1$ be two disjoint subsets of $\Z^d$ such that
%$\inf\{\max_{1\le k\le d} |i_k-j_k|:\, i \in \Gamma_0, j \in \Gamma_1\}=1$.
%\sidecomment{A stronger assumption is that $\{\Gamma_0,\Gamma_1\}$
%is a partition of $\Z^d$. But, perhaps, the weaker assumption might be
%enough.}
%Then there are constants $R_1>0$ and $0<R_2 < 1$ such that
Then, there are constants $D_1>0$ and $0<D_2 < 1$ such that
 \begin{align*}
\nu \|\nabla_{\Gamma_i} (\E^{\Gamma_j} f)\|^2 
& \leq D_1 \nu \| \nabla_{\Gamma_i} f\|^2 + D_2 \nu \| \nabla_{\Gamma_j} f\|^2.
\end{align*}
for $i\neq $j, $i,j-0,1$
\end{lemma}
\begin{proof}
Fix a site $i \in \Gamma_1$. If then we denote by $\partial\{i\}$ the set  
of the   neighbors of $i$,  since $\partial \{i\} \subset \Gamma_0$,
we get   $\E^{\Gamma_0} f = \E^{\Gamma_0\setminus \partial\{i\}} \E^{\partial\{i\}}f$.
Hence  for a vector field $X^\alpha_i$ being   one of the H\"ormander generators
 of the spin space $\Spin^{\{i\}}$, we can write
$X^\alpha_i \E^{\Gamma_0} f = \E^{\Gamma_0\setminus \partial\{i\}} X^\alpha_i(\E^{\partial\{i\}}f)$. For $i \in \Gamma_1$, we have
\[
\nu\|\nabla_{i} (\E^{\Gamma_0} f)\|^2=\sum_{a=1}^{n}\nu(X^\alpha_{i} (\E^{\Gamma_0} f))^2
=\sum_{a=1}^{n}\nu(\E^{\Gamma_0\setminus \partial\{i\}} X^\alpha_i(\E^{\partial\{i\}}f) )^2
.\]
By Jensen's inequality
\[
\nu\|\nabla_{i} (\E^{\Gamma_0} f)\|^2
\leq \sum_{a=1}^{n}\nu\E^{\Gamma_0\setminus \partial\{i\}}( X^\alpha_i(\E^{\partial\{i\}}f) )^2= \sum_{a=1}^{n}\nu( X^\alpha_i(\E^{\partial\{i\}}f) )^2=\nu \|\nabla_i(\E^{\partial\{i\}}f)\|^2
\]
since  $\nu \E^{\Gamma_0 \setminus \partial\{i\}}=\nu$. Finally, 
summing over $i \in \Gamma_1$ we have
\[
\nu\|\nabla_{\Gamma_1} (\E^{\Gamma_0} f)\|^2
\le \sum_{i\in \Gamma_1} \nu \|\nabla_i(\E^{\partial\{i\}}f)\|^2
.\]
If we denote  $\partial\{i\} 
=\{j_1, j_2, j_3, j_4\}$ the neighbors of $i$ and use  Lemma \ref{sweep1}  we have
\[
\nu\|\nabla_i (\E^{\partial\{i\}}f)\|^2
\le Q_1 \nu\|\nabla_i \E^{\{j_2,j_3,j_4\}} f\|^2
+ Q_2 \nu\|\nabla_{j_1} \E^{\{j_2,j_3,j_4\}} f\|^2
.\]
Since the neighbors of $i$ all have distance bigger than one from each other, 
$\nabla_{j_1} \E^{\{j_2,j_3,j_4\}} f
= \E^{\{j_2,j_3,j_4\}} \nabla_{j_1} f$ and so, for the second term, by Jensen's inequality
\[
\nu\|\nabla_{j_1} \E^{\{j_2,j_3,j_4\}} f\|^2
\le \nu  \|\nabla_{j_1} f\|^2.
\]
For the first term we can use once more Lemma \ref{sweep1},
\[
\nu\|\nabla_i \E^{\{j_2,j_3,j_4\}} f\|^2
\le Q_1 \nu\|\nabla_i \E^{\{j_3,j_4\}} f\|^2
+ Q_2 \nu \|\nabla_{j_2} \E^{\{j_3,j_4\}} f\|^2.
\]
If we continue in this manner, we finally  obtain 
\begin{align*}
\nu\|\nabla_i (\E^{\partial\{i\}}f)\|^2
\le Q_1^4 \nu\|\nabla_i f\|^2 + Q_1^3 Q_2 \sum_{j \in \partial\{i\}}
\|\nabla_j f\|^2.
\end{align*}
and now  summing  over  $i \in \Gamma_1$,
\[
\nu \|\nabla_{\Gamma_1}(\E^{\Gamma_0} f)\|^2
\le Q_1^4 \nu\|\nabla_{\Gamma_0} f\|^2 
+ 4 Q_1^3 Q_2 \nu\|\nabla_{\Gamma_1} f\|^2.
\]
Substituting the
 expressions for $Q_1$ and $Q_2$ from Lemma \ref{sweep1} (recall that $Q_2<1$ for $J$ small enough) we obtain the inequality with $D_1 :=Q_1^4$
and $D_2 := 4 Q_1^3 Q_2<1$, provided that
$J < (80 \lambda (c+2cG_0+2G_0))^{-1/2}$.
\end{proof}
To continue we need to define the symbol $\D$ as follows: define $\D^0 f=f$. For   $n$ odd $\D^n := \E^{\Gamma_0} \D^{n-1}$, while for $n$ even   $\D^n := \E^{\Gamma_1} \D^{n-1}$. 
In that way   $\D^{n}$ for $n$  even   is a functional that does not depend on $x_{\Gamma_1}$, while     $\D^{n} f$ for $n$ odd does not depend on $x_{\Gamma_0}$. 

\begin{lemma}\label{conv}
Assuming (\ref{Hphi})-(\ref{HgradV2}),  we  have
that $\lim_{n \to \infty} \D^n f = \nu f$, $\nu$-a.e.
\end{lemma}
\begin{proof}

For $k=0,1$, we compute 
\[\nu (\D^{n} f-  \D^{n+1} f)^2=\nu \left(\E^{\Gamma_k}(\D^{n} f-\E^{\Gamma_k} \D^{n} f)^2\right) \le c\,\E^{\Gamma_0}\|\nabla_{\Gamma_k}\D^{n} f\|^2\]
where above we used the  spectral gap inequality for the product measure  $\E^{\Gamma_k},k=0,1$,  which follows
from the product property of the spectral gap inequality, which itself follows from the log-Sobolev inequality.  
If  we use $n$ times Lemma \ref{sweepGam1} to estimate the last term we obtain
\[
\nu[(\D^n f-\D^{n+1} f)^2] \le  (Q_1 +Q_{2})Q_2^{n-1} \nu \|\nabla   f\|^2 \le T^n,
\]
for some $0<T<1$ that depends on $f$.
From the last inequality the lemma follows directly   using relevant  techniques developed in \cite{G-Z}. For a detailed exposition one can also look on a similar application in \cite{Pa1} and \cite{Pa2}. 
\end{proof}

\section{Second set of sweeping-out inequalities}
In this section we show further sweeping out inequalities.\begin{lemma}\label{sweep2}
Assume (\ref{Hphi})-(\ref{HgradV2}) and that the log-Sobolev inequality holds uniformly  for $\E^{i,\omega}$. Then for  $i \sim j$,  there exist constants  $R_1>0$ and $0<R_2<1$, so that
\[
\nu\| \nabla_i \sqrt{\E^j f^2}\|^2 
\leq R_1 \nu \|\nabla_i f\|^2 + R_2 \nu\|\nabla_j f\|^2
+ R_2 \sum_{k \sim j,\, k \neq i} \nu \|\nabla_k f\|^2
.\]
\end{lemma}
\begin{proof}
For two neighboring sites $i$ and $ j$, compute
\begin{equation}
\label{above} \| \nabla_i \sqrt{\E^j f^2}\|^2 = \sum_{\alpha=1}^n
(X_i^\alpha \sqrt{\E^j f^2})^2=
\sum_{\alpha=1}^n \frac{(X_i^\alpha \E^j f^2)^2}{4\E^j f^2}.
\end{equation}
 The numerator can be  estimated as in (\ref{F1a}), 
\[
(X_i^\alpha \E^j f^2)^2
\le 2 (\E^j(X_i f^2))^2 + 2(\int(X_i^\alpha \theta_j) f^2 dx_j)^2,\]
where we recall that $\theta_j$ stands for the density of $\E^j$. 
Then, if we apply Leibnitz' rule, Cauchy-Schwarz   as well as  Jensen inequalities to bound  
 the first term, while we   bound the second summand  from (\ref{F1b}) we obtain
\[
(X_i^\alpha \E^j f^2)^2  \,
\le 4 (\E^j f^2) \E^j(X_i^\alpha f)^2
+ 2 J_{ji}^2 \cov_{\E^j}[f^2,\, X_i^\alpha V(x_j, x_i)] ^2,
\]
where we denoted $\cov_\mu(f, g) = \mu(fg)-\mu(f)\mu(g) $, for any probability measure $\mu$. From this together with  (\ref{above})   we have
\[
\|\nabla_i \sqrt{\E^j f^2}\|^2  \le \E^j\|\nabla_i f\|^2 + \frac{J^2}{2} 
\sum_\alpha \frac{\cov_{\E^j}[f^2,\, X_i^\alpha V(x_j, x_i)] ^2}{\E^j f^2}.
\]
We can bound the second term on the right hand side by the following estimation 
\[
(\cov_\mu(f^2, g))^2 \le 8\, (\mu f^2) \mu [(f-\mu f)^2 (g^2 + \mu g^2)]
\] 
(see \cite{Pa1} for a proof). We then  get
\begin{equation}
\label{starting}
\|\nabla_i \sqrt{\E^j f^2}\|^2  \le \E^j\|\nabla_i f\|^2 + 4J^2  
\E^j\big \{(f-\E^j f)^2 (\|\nabla_i V\|^2 + \E^j \|\nabla_i V\|^2) \big\}.
\end{equation}
The following inequality presented in  \cite{D-S} provides a bound for the second term:

Let $Z, W$ be measurable  functions with $Z>0$.
Under appropriate integrability assumptions
\[
\mu(DB) \le \mu(B \log(B/ \mu B)) + (\mu B) \log \mu(e^D).
\]
In our case we  apply the inequality for  $D=\epsilon^{-1} (f-\E^j f)^2$,
and $B = \epsilon(\|\nabla_i V\|^2 + \E^j \|\nabla_i V\|^2)$, so that
\begin{align}\nonumber
\E^j\big \{(f-\E^j f)^2 (\|\nabla_i V\|^2 + \E^j \|\nabla_i V\|^2) \big\}
\leq &
\epsilon^{-1} \E^j \big[(f-\E^j f)^2 \log \frac{(f-\E^j f)^2}{\E^j (f-\E^j f)^2}\big]  \\
+ &\epsilon^{-1} \underbrace{(\E^j (f-\E^j f)^2) 
\log \E^j e^{\epsilon (\|\nabla_i V\|^2 +\E^j \|\nabla_i V\|^2)}}_{\O(i,j)}.
\label{useentropic}
\end{align}
The first term on the right is bounded  by the log-Sobolev
inequality 
\begin{equation}
\label{useLS}
\E^j \big[(f-\E^j f)^2 \log \frac{(f-\E^j f)^2}{\E^j (f-\E^j f)^2}\big]
\le c \E^j\|\nabla_j (f-\E^j f)\|^2= c \E^j\|\nabla_j f\|^2
\end{equation}
while  for the second term, we   use
%\begin{align}
%\|\nabla_i \sqrt{\E^j f^2}\|^2  
%&\le \E^j\|\nabla_i f\|^2 + \frac{4J^2}{\epsilon} 
%\left\{\E^j \big[\hf^2 \log (\hf^2/\E^j \hf^2)\big]
%+ (\E^j \hf^2) 
%\log \E^j e^{\epsilon (\|\nabla_i V\|^2 +\E^j \|\nabla_i V\|^2)}\right\}
%\nonumber
%\\
%&\le \E^j\|\nabla_i f\|^2 
%+ \frac{4 c J^2}{\epsilon} \E^j \|\nabla_j f\|^2
%+ \frac{4J^2}{\epsilon} (\E^j \hf^2)
%\big\{
%\underbrace{\epsilon \E^j \|\nabla_i V\|^2 
%+ \log \E^j e^{\epsilon \|\nabla_i V\|^2}}_{:=T}
%\big\}.
%\label{TTT}
%\end{align}
the two bounds  (\ref{HgradV1}) and (\ref{HgradV2}) for the interactions, so that, 
\begin{align*}
\O(i,j)=&\E^j (f-\E^j f)^2 \big\{\epsilon \E^j \|\nabla_i V\|^2 
+ \log \E^j e^{\epsilon \|\nabla_i V\|^2}\big\}  
\\ \leq &(\epsilon+1)\lambda \E^j (f-\E^j f)^2 +\epsilon \lambda (\E^j (f-\E^j f)^2) \E^j[\dd(x_j)^r]
+ \\ &  +(\epsilon+1) \lambda \E^j [(f-\E^j f)^2 \dd(x_i)^r]  
+ \lambda \sum_{\ell \sim j, \ell \neq i}  \E^j[(f-\E^j f)^2 \dd(x_\ell)^r].
\end{align*}
Taking expectations with  respect to $\nu$ 
\begin{multline*}
\nu (\O(i,j))
\le \nu\bigg\{
(\epsilon+1)\lambda \E^j (f-\E^j f)^2 + \epsilon \lambda (\E^j (f-\E^j f)^2) \E^j[\dd(x_j)^r]
+ \\ (\epsilon+1) \lambda \E^i [(f-\E^j f)^2 \dd(x_i)^r]  
+ \lambda \sum_{\ell \sim j, \ell \neq i}  \E^\ell[(f-\E^j f)^2 \dd(x_\ell)^r]
\bigg\}.
\end{multline*}
The variance  $\E^j (f-\E^j f)^2$ appearing on the first term is upper bounded by the spectral gap inequality
 $\E^j (f-\E^j f)^2 \le c \E^j \|\nabla_j f\|^2$, while 
 for the remaining three  terms involving $\dd$ we can use  the coercive inequality obtained in  Lemma \ref{coers}. Then, the first of these terms is just 
$\E^j[\dd(x_j)^r] \le G_0$, while the second one becomes 
\[
\E^i [(f-\E^j f)^2 \dd(x_i)^r]
\le G_0 \E^i \|\nabla_i (f-\E^j f)\|^2 + G_0 \E^i (f-\E^j f)^2
\leq (c+1) G_0 \E^i \|\nabla_i f\|^2,
\]
where above we also used again the spectral gap inequality for $\E^i$. In the same way, we compute 
\[\E^\ell [(f-\E^j f)^2 \dd(x_\ell)^r]
\le (c+1) G_0 \E^\ell \|\nabla_\ell f\|^2.\]
Collecting all the above together,
\[
\nu (\O(i,j)) \le \lambda (\epsilon+1+\epsilon G_0) \nu\|\nabla_j f\|^2
+ \lambda (\epsilon+1)(c+1)G_0 \nu\|\nabla_i f\|^2
+ \lambda(c+1)G_0 \sum_{\ell \sim j, \ell \neq i} \nu\|\nabla_\ell f\|^2.
\]
Combining  the last one together with \eqref{starting}-\eqref{useLS} proves the lemma with constants
$R_1 = 1+ 4J^2 \lambda(1+\epsilon^{-1})(c+1)G_0$
and $R_2 = 4J^2 \max\{\epsilon^{-1} c+\lambda(1+\epsilon^{-1} +G_0),\,
\lambda \epsilon^{-1} (c+1)G_0\}<1$, provided that 
 $J < 1/2\max\{\epsilon^{-1} c+\lambda (1+\epsilon^{-1} +G_0),\,
\lambda \epsilon^{-1} (c+1)G_0\}^{1/2}$.
\end{proof}

\begin{lemma}\label{sweepGam2}
Assume conditions     (\ref{Hphi})-(\ref{HgradV2})   and that  the log-Sobolev inequality holds uniformly for $\E^{i,\omega}$. There are constants $G_1>0$ and $0<G_2 < 1$ such that
 \begin{align*}
\nu\|\nabla_{\Gamma_i} \sqrt{\E^{\Gamma_j} f^2}\|^2
&\le G_1\, \nu\|\nabla_{\Gamma_i} f\|^2 + G_2\, \nu\|\nabla_{\Gamma_j} f\|^2
.\end{align*}for $i\neq j, i,j=0,1$
\end{lemma}
\begin{proof}
Consider any two subsets of $\Z^2$, $A$ and $B$  at lattice distance two or bigger from each other, and any $i \in \Z^2$ such that 
$\partial\{i\} \cap A = \varnothing$. Then, for   a gradient   $\nabla_i = (X^\alpha_i, . . . , X^\alpha_n)$, after using twice the  Leibnitz' rule we can 
  write
\begin{align*}\nu\|\nabla_i \sqrt{\E^{A\cup B} f^2}\|^2=\nu \sum_{a=1}^n \left(
X^\alpha_i \sqrt{\E^{A\cup B} f} \right)^2
=\nu \sum_{a=1}^n \left(\frac{ \E^A [ \sqrt{\E^B f} X^\alpha_i \sqrt{\E^B f}]}{\sqrt{\E^{A\cup B}f}}\right)^2
.\end{align*}
By Cauchy-Schwarz inequality 
\begin{align}\label{QS}\nu\|\nabla_i \sqrt{\E^{A\cup B} f^2}\|^2
 \leq \nu\sum_{a=1}^n \frac{  (\E^A \E^B f) \, \E^A(X^\alpha_i \sqrt{\E^B f})^2}{\E^{A\cup B}f}=\nu \|\nabla_i \sqrt{\E^{B} f^2}\|^2.
\end{align}
Using (\ref{QS}) we can  bound the left side of the statement:
\begin{equation}
\label{backtoit}
\int \|\nabla_{\Gamma_1}  \sqrt{\E^{\Gamma_0} f^2}\|^2 d\nu
\le \sum_{i\in\Gamma_1} 
\underbrace{\int \|\nabla_i  \sqrt{\E^{\partial\{i\}}f^2}\|^2 d\nu}_{:=A_1(i)}
.\end{equation}
For  $i\in \Z^2$, we will upper bound $A_i(i)$. 

Consider a node $i \in \Gamma_1$ and denote its four neighbors $\{\sim i \}$ by $i_1, i_2,i_3,i_4$. Then denote $i_{mn}$ the node at distance two from $i$ that neighbours with $i_m$ and $i_n$ and as $i_m'$ the node at distance two from $i$ that neighbours only with $i_m$ from $\{\sim i\}$.  
We begin with the sweeping out inequality from   Lemma \ref{sweep2}. We have
\begin{align}
\label{T1}\nonumber
A_1(i):&
\le R_1 \underbrace{\int \|\nabla_i \sqrt{\E^{\{i_2,i_3,i_4\} } f^2}\|^2 d\nu}_{:=A_2(i)}+
\\ &+ R_2\int \|\nabla_{i_1}  \sqrt{\E^{\{i_2,i_3,i_4\} } f^2}\|^2 d\nu
+ R_2 \sum_{k\in\{i_1',i_{12},i_{14}\}} \int \|\nabla_k  \sqrt{\E^I f^2}\|^2 d\nu
.\end{align}
In order to pass the gradient inside the expectation, we will apply inequality   (\ref{QS})  three times in the second term on the right hand side above. That will give 
\[
\int \|\nabla_{i_1} \sqrt{\E^{\{i_2,i_3,i_4\} }f^2}\|^2 d\nu \leq \int \|\nabla_{i_1} f\|^2 d\nu
.\]%In the third term, $k$ ranges over $\{i_{12}, i_{14}, i_1'\}$.
Concerning the      last sum in (\ref{T1})  we estimate each one of the three terms separately, as shown below. 
Using again inequality  (\ref{QS})    and  Lemma \ref{sweep2} we get
\begin{align*}
\int \|\nabla_{i_{12 }}  \sqrt{\E^{\{i_2,i_3,i_4\} } f^2}\|^2 d\nu
\le&\int \|\nabla_{i_{12 }}  \sqrt{\E^{i_2} f^2}\|^2 d\nu
\le
R_1\int \|\nabla_{i_{12}}f\|^2 d\nu
+R_2\int \|\nabla_{i_2}f\|^2 d\nu+  \\  &
+R_2\sum_{\ell \sim i_2, \ell \neq i_{12}} \int \|\nabla_\ell f\|^2 d\nu
.\end{align*}
For the other term we have
\begin{align*}
\int \|\nabla_{i_{14 }}  \sqrt{\E^{\{i_2,i_3,i_4\} } f^2}\|^2 d\nu
\le & \int \|\nabla_{i_{14 }}  \sqrt{\E^{i_4} f^2}\|^2 d\nu
\le
R_1\int \|\nabla_{i_{14}}f\|^2 d\nu
+R_2\int \|\nabla_{i_4}f\|^2 d\nu + \\  &
+R_2\sum_{\ell \sim i_4, \ell \neq i_{14}} \int \|\nabla_\ell f\|^2 d\nu
\end{align*}
while using  (\ref{QS}) once more
\[
\int \|\nabla_{i_1'}  \sqrt{\E^{\{i_2,i_3,i_4\} } f^2}\|^2 d\nu
=\int \|\nabla_{i_1'}  \sqrt{\E^{\{i_2,i_3,i_4\} } f^2}\|^2 d\nu
\le \int \|\nabla_{i_1'} f\|^2 d\nu
.\]
We gather together the last four displays  to  bound at first (\ref{T1}) and then 
   sum over $i \in \Gamma_1$ to upper bound (\ref{backtoit}). Then we obtain
\[
 \sum_{i \in \Gamma_1} A_1(i) 
\le R_1 \sum_{i \in \Gamma_1} A_2(i)
+ R_2 (2 G_2+1) \int \|\nabla_{\Gamma_0}f\|^2 d\nu
+ R_2 (6G_2 + 2G_1 + 1) \int \|\nabla_{\Gamma_1}f\|^2 d\nu
.\]
%\begin{multline*}
%\sum_{i \in \Gamma_1} T_1(i)
%\le
%G_1 \sum_{i \in \Gamma_1} T_2(i)
%+ G_2 \sum \III{\nabla_{i_1} f}^2
%+ G_1 G_2 \sum (\III{\nabla_{i_{12}}f^2}^2+\III{\nabla_{i_{14}}f^2}^2)
%\\
%+ G_2^2 \sum(\III{\nabla_{i_2}f^2}^2 + \III{\nabla_{i_4}f^2}^2)
%+ G_2 \sum \III{\nabla_{i_1'} f}^2
%+ G_2^2 \sum(\sum_{\ell \sim i_2, \ell \neq i_{12}} \III{\nabla_\ell f^2}^2
%+ \sum_{\ell \sim i_4, \ell \neq i_{14}} \III{\nabla_\ell f^2}^2)
%\\
%= G_1 \sum_{i \in \Gamma_1} T_2(i)
%+ G_2 S(\Gamma_0)
%+2 G_1G_2 S(\Gamma_1)
%+ 2 G_2^2 S(\Gamma_0)
%+ G_2 S(\Gamma_1)
%+ 6 G_2^2 S(\Gamma_1)
%\\
%= G_1 \sum_{i \in \Gamma_1} T_2(i)
%+ (G_2+2 G_2^2) S(\Gamma_0) + (2G_1G_2 + G_2 + 6G_2^2) S(\Gamma_1),
%\end{multline*}
%where $S(\Gamma_i) :=\III{\nabla_{\Gamma_i} f}^2$, $i=0,1$.
In order to estimate $A_2(i)$, we proceed in a similar way. We can write
\begin{align}\nonumber
\label{T2}
 A_2(i)
= & \int \|\nabla_i \sqrt{\E^{i_2} \E^{\{i_3,i_4\}} f^2}\|^2 d\nu
\le 
R_1 \underbrace{\int \|\nabla_i  \sqrt{\E^{\{i_3,i_4\}} f^2}\|^2 d\nu}_{:=A_3(i)}
+ R_2\int \|\nabla_{i_2}  \sqrt{\E^{\{i_3,i_4\}} f^2}\|^2 d\nu + \\ &
+ R_2 \sum_{k\in\{i_2', i_{12}, i_{23}\}} \int \|\nabla_k  \sqrt{\E^{\{i_3,i_4\}} f^2}\|^2 d\nu
.\end{align}
%We let $T_3(i) := \III{\nabla_i \widehat\E^J f^2}^2$ and use (SQ) for the
As before we use (\ref{QS})  for the second term,
$
\int \|\nabla_{i_2}  \sqrt{\E^{\{i_3,i_4\}}f^2}\|^2 d\nu
\le\int \|\nabla_{i_2}f\|^2 d\nu
$,
while we will  
estimate   the   terms of the last sum of (\ref{T2}) separately as shown below:
\begin{align*}
\int \|\nabla_{i_{12}}  \sqrt{\E^{\{i_3,i_4\}} f^2}\|^2 d\nu \le &
\int \|\nabla_{i_{12}}f\|^2 d\nu, \  \  \   \int \|\nabla_{i_2'}  \sqrt{\E^{\{i_3,i_4\}} f^2}\|^2 d\nu \le \int \|\nabla_{i_2'} f\|^2 d\nu
\\
\int \|\nabla_{i_{23}}  \sqrt{\E^{\{i_3,i_4\}} f^2}\|^2 d\nu \le &
\int \|\nabla_{i_{23}} \sqrt{\E^{i_3} f^2}\|^2 d\nu
\le R_1\int \|\nabla_{i_{23}}f\|^2 d\nu
+ R_2\int \|\nabla_{i_3}f\|^2 d\nu + \\  &
+ R_2 \sum_{\ell\sim i_3,\ell \neq i_{23}}\int \|\nabla_\ell f\|^2 d\nu
\end{align*}
where for the second we applied Lemma \ref{sweep2}. Substituting these into (\ref{T2}) and then  taking the sum   over $i \in \Gamma_1$
\[
 \sum_{i\in \Gamma_1} A_2(i)
\le R_1 \sum_{i\in \Gamma_1} A_3(i) 
+ R_2(R_2+1) \int \|\nabla_{\Gamma_0}f\|^2 d\nu
+ R_2(4R_2+R_1R_2+1) \int \|\nabla_{\Gamma_1}f\|^2 d\nu
.\]
%%%OOOOOOOOOOOOOOOOOOOOOOOOOOOOOOOOOOOOOOOOOOOOOOOOOOOOOOOOOOOOOOOOOOOOOOOOOOOOOOOOOOOOOOOOOOOOOOOOOOOOOO
%\begin{multline*}
%\sum_{i\in \Gamma_1} T_2(i) \le G_1 \sum T_3(i) 
%+ G_2 \sum\III{\nabla_{i_2}f}^2
%+G_2 \sum \III{\nabla_{i_{12}}f}^2
%+ G_1 G_2 \sum \III{\nabla_{i_{23}}f}^2
%+ G_2^2 \III{\nabla_{i_3}f}^2
%\\
%+ G_2^2 \sum \sum_{\ell\sim i_3,\ell \neq i_{23}}\III{\nabla_\ell f}^2
%+ G_2 \III{\nabla_{i_2'} f}^2
%\\
%= G_1 \sum T_3(i) + G_2 S(\Gamma_0) + G_2 S(\Gamma_1)
%+ G_1 G_2 S(\Gamma_1) + G_2^2 S(\Gamma_0)
%+ 3 G_2^2 S(\Gamma_1) + G_2 S(\Gamma_1)
%\\
%= G_1 \sum T_3(i) + (G_2+G_2^2) S(\Gamma_0) 
%+ (G_2G_2+3G_2^2+2G_2) S(\Gamma_1).
%\end{multline*}
Finally, we calculate $A_3(i)$   similarly.
\begin{align*}
 A_3(i)
=&\int \|\nabla_i  \sqrt{\E^{i_3} \E^{i_4} f^2}\|^2 d\nu
\le R_1\int \|\nabla_i  \sqrt{\E^{i_4}f^2}\|^2 d\nu
+ R_2 \int \|\nabla_{i_3}  \sqrt{\E^{i_4}f^2}\|^2 d\nu +\\ &
+ R_2 \sum_{k\in\{i_3',i_{23},i_{24}\}}\int \|\nabla_k \sqrt{\E^{i_4}f^2}\|^2 d\nu
\end{align*}
with the remaining terms estimated as shown below with the use of (\ref{QS}) and Lemma  \ref{sweep2}.
\begin{align*}&
\int \|\nabla_i  \sqrt{\E^{i_4}f^2}\|^2 d\nu \le
R_1 \int \|\nabla_{i} f\|^2 d\nu
+ R_2\int \|\nabla_{i_4} f\|^2 d\nu
+ R_2 \sum_{k \sim i_4, k \neq i} \int \|\nabla_{k} f\|^2 d\nu
\\   &
\int \|\nabla_{i_3} \sqrt{\E^{i_4}f^2}\|^2 d\nu \le \int \|\nabla_{i_3} f\|^2 d\nu 
,   \   \   \
\int \|\nabla_{i_{23}} \sqrt{\E^{i_4}f^2}\|^2 d\nu \le \int \|\nabla_{i_{23}}f\|^2 d\nu
\\   &
\int \|\nabla_{i_{34}} \sqrt{\E^{i_4}f^2}\|^2 d\nu \le
R_1\int \|\nabla_{i_{34}} f\|^2 d\nu 
+ R_2\int \|\nabla_{i_4} f\|^2 d\nu
+ R_2 \sum_{\ell\sim i_4, \ell\neq i_{34}}\int \|\nabla_{\ell} f\|^2 d\nu
.\end{align*}
This leads to\[
 \sum_{i\in \Gamma_1} A_3(i)
%\le R_1 \sum_{i\in \Gamma_1} A_1(i)
\le  R_2(R_2+R_1+1)\int \|\nabla_{\Gamma_0}f\|^2 d\nu
+ (3R_2^2+4R_1R_2+R_2+R_1^2) \int \|\nabla_{\Gamma_1}f\|^2 d\nu
.\]
%\[
%{\red \sum_{i\in \Gamma_1} T_2(i)}
%\le G_1 \sum_{i\in \Gamma_1} T_3(i) 
%+ G_2(G_2+1) \III{\nabla_{\Gamma_0}f}^2
%+ G_2(4G_2+G_1G_2+1) \III{\nabla_{\Gamma_1}f}^2.
%\]
%\[
%{\red \sum_{i \in \Gamma_1} T_1(i)} 
%\le G_1 \sum_{i \in \Gamma_1} T_2(i)
%+ G_2 (2 G_2+1) \III{\nabla_{\Gamma_0}f}^2
%+ G_2 (6G_2 + 2G_1 + 1) \III{\nabla_{\Gamma_1}f}^2.
%\]
To finish the proof,  we substitute the terms involving the sums of $A_1(i),A_2(i)$ and $A_3(i)$ to one another and then
back to (\ref{backtoit}) to obtain  the second inequality in
the statement with constants
$G_1 = 4{R_1}^{2}{R_2}^{2}+4{R_1}^{3}R_2+{R_1}^{2}R_2+
{R_1}^{4}+4R_1{R_2}^{2}+3R_1G_2+6{R_2}^{2}+R_2$ and  
$G_2 = R_2 ( {R_1}^{2}R_2+{R_1}^{3}+{R_1}^{2}+R_1R_2+R_1+2\,R_2+1).$ 
Since constant $R_2 = 4J^2 \max\{\epsilon^{-1} c+\lambda (1+\epsilon^{-1} +G_0),\,
\lambda \epsilon^{-1} (c+1)G_0\}$,  we can choose $J$ sufficiently small such that $R_2$ is small enough so that $G_2<1$.

\end{proof}

%Therefore $\ent_{\E^\Lambda}(g)$
%depends on $x \in \Spin^{\Z^d}$ only through $x_{\Lambda^c}$.
%Integrating the last display with respect to the Gibbs measure $\nu$
%and using $\nu \E^\Lambda =\nu$ we obtain
%\begin{equation}
%\label{recit}
%\nu[ g \log g] = \nu[\ent_{\E^\Lambda}(g)] + 
%\nu[(\E^\Lambda g) (\log \E^\Lambda g)].
%\end{equation}

\section{The log-Sobolev inequality for the Gibbs measure.}

Now we can  prove the main result stated in Theorem \ref{thmGENERAL}.
 \begin{proof}

 Denote $\ent_{\mu}(  h):=\mu \left(h\log \frac{h}{\mu h}\right)$ the entropy of $h$ with respect to the measure $\mu$. At the beginning we need the following equality
\begin{align}\nonumber
\label{Qbase}
\D^n[h \log h] 
= &\sum_{m=0,\ m\ \text{even} }^{n-1} \D^{n-m}  [\ent_{\E^{\Gamma_0}}(\D^m h)]
+ \sum_{m=0,\  m\ \text{odd}}^{n-1} \D^{n-m-1} [\ent_{\E^{\Gamma_1}}( \D^m h)]
+ \\  &+ (\D^n h)\log(\D^n h)
\end{align}
for  all $n \geq 0$.
The statement can be easily proven by induction (see \cite{Pa2}).
From   Lemma \ref{conv}, we obtain the limits  $\D^n[\lambda(h)] \to \nu[\lambda(h)]$
and $\lambda^n(\D^n h) \to \nu[h]$, $\nu$-a.e. From these two limits  together with 
  Fatou's lemma, (\ref{Qbase})  becomes
\begin{align}
\ent_\nu(h) &\le \liminf_{n \to \infty}\bigg\{
\sum_{^{m=0,\ m\ \text{even}} }^{n-1} \nu[\ent_{\E^{\Gamma_0}}(\D^m h)]
+ \sum_{m=0,\ m\ \text{odd}}^{n-1} \nu[\ent_{\E^{\Gamma_1}}(\E^{\Gamma_0}\D^m h)]
\label{entnug}
\bigg\}
\end{align}
where in order to obtain the last equality we  also used that $\nu(\D^{n-m-1} f)=\nu (f)$ since $\nu$ is a Gibbs  measure.
Now substitute $h=f^2$ and observe that for the first term of the first sum we have
\[
\nu[\ent_{\E^{\Gamma_0}}(f^2)] \le c \nu \|\nabla_{\Gamma_0} f\|^2,
\]
the reason being that the log-Sobolev inequality holds for $\E^{\Gamma_0}$,
since the elements in $\Gamma_0$ are at lattice distance strictly greater
 than one from one another, while  $\nu \E^{\Gamma_0} = \nu$. Similarly, for $m \ge 1$, we apply
first the log-Sobolev for $\E^{\Gamma_0}$ and $\E^{\Gamma_1}$ and then
integrate with respect to $\nu$ ($i\neq j$, $i,j=0,1$):
\begin{align*}
\nu[\ent_{\E^{\Gamma_i}}(\D^{m} f^2)]
&\le c \nu \|\nabla_{\Gamma_i} \sqrt{\D^m f^2}\|^2
\le c [G_1 C_2^{2m-1}  \nu \|\nabla_{\Gamma_j} f\|^2 
+G_2^{2m} \nu \|\nabla_{\Gamma_i} f\|^2 ]
,
\end{align*}
for $m$ odd and even respectively,  where,  for the last inequalities we used inductively $m$ times Lemma \ref{sweepGam2}.
Finally, substituting the last two  in (\ref{entnug}), we obtain that
\[
\ent_\nu(f^2) \le 
\frac{c (1+G_2)}{1-G_2^2} \nu \|\nabla_{\Gamma_1} f\|^2 
+ \frac{c (1+G_2)}{1-G_2^2} \nu \|\nabla_{\Gamma_0} f\|^2 
\]
since   $0< G_2 < 1$.
This  proves the log-Sobolev inequality for $\nu$. 
\end{proof}

\bibliographystyle{alpha}

\end{document}